\documentclass[]{siamart171218}

\usepackage{blindtext}
\usepackage{multicol}
\usepackage{amsmath}

\DeclareMathOperator*{\argmin}{arg\,min}
\usepackage{amsfonts}
\usepackage{graphicx}
\usepackage{hyperref}
\usepackage{amsmath,amssymb,setspace}
\usepackage{algorithm}
\usepackage{graphics}
\usepackage{booktabs}
\usepackage{algpseudocode}
\usepackage{tikz}
\usepackage{tikz-3dplot}
\usepackage{bm}
\newcommand{\bs}{\boldsymbol}

\newcommand{\mb}{\mathbb}
\newcommand{\mf}{\mathbf}
\newcommand{\Depth}{2}
\newcommand{\Height}{2}
\newcommand{\Width}{2}

\usepackage{color,amsmath}

\headers{Covid-19 Analysis Using Tensor Methods}{D. Dulal, R. Goudarzi Karim, and C. Navasca}

\title{Covid-19 Pandemic Data Analysis using Tensor Methods \thanks{Submitted to the editors December 29, 2022.
\funding{This work was funded by National Science Foundation under Grant No. DMS-1439786 and Grant No. MCB-2126374.}}}

\author{Dipak Dulal\thanks{Department of Mathematics, University of Alabama at Birmingham, Birmingham, AL 35294
  (\email{dpdulal@uab.edu)}.}
\and Ramin Goudarzi Karim\thanks{Department of Computational and Information Sciences, Stillman College, Tuscaloosa, AL 35401
  (\email{rkarim@stillman.edu)}.}
\and Carmeliza Navasca\thanks{Department of  Mathematics, University of Alabama at Birmingham, Birmingham, AL 35294
  (\email{cnavasca@uab.edu)}.}}

\begin{document}

\maketitle

\begin{abstract}
In this paper, we use tensor models to analyze Covid-19 pandemic data. First, we use tensor models, canonical polyadic and higher-order Tucker decompositions, to extract patterns over multiple modes. Second, we implement a tensor completion algorithm using canonical polyadic tensor decomposition to predict spatiotemporal data from multiple spatial sources and to identify Covid-19 hotspots. We apply a regularized iterative tensor completion technique with a practical regularization parameter estimator to predict the spread of Covid-19 cases and to find and identify hotspots. Our method can predict weekly and quarterly Covid-19 spreads with high accuracy. Third, we analyze Covid-19 data in the US using a novel sampling method for alternating least-squares. Moreover, we compare the algorithms with standard tensor decompositions it terms of their interpretability, visualization and cost analysis. Finally, we demonstrate the efficacy of the methods by applying the techniques to New Jersey’s Covid-19 case tensor data.

\end{abstract}


\section{Introduction}
Tensor decomposition is a powerful tool in data analysis, computer vision, scientific computing, machine learning and many other fields. In fact, tensor models for dimensionality reduction have been highly effective in machine learning applications like classification and regression. Tensor decomposition on its own have been successful in many modern big data analysis \cite{karim2020accurate,fat_nav,fatou,nav_wang}.

In this work, we focus on analyzing Covid-19 pandemic data \cite{CovidData,PopData} by using tensor decomposition models. Our goals are to predict future infection, locate and identify hotspots from data-set coming from multiple sources across different spatial regions over time. One is interested in determining where and when changes occur in the pattern. The strategy is set up an optimization which separates the data in following: let $\mathcal{Y} = \mathcal{L} + \mathcal{S}$ where $\mathcal{Y}$ is the given tensor,  $\mathcal{L}$ is a low rank reconstructed tensor of $\mathcal{Y}$ and $\mathcal{S}$ is the sparse tensor. In video processing, the original video is separated into its background and foreground subspace to detect anomalous activities \cite{karim2020accurate,WangNa2}. The tensor $\mathcal{L}$ is the background, and $\mathcal{S}$ is the foreground. The sparse tensor $\mathcal{S}$ can provide anomalous activities. Similarly, $\mathcal{S}$ will contain hotspot occurrence. To achieve this separation, we formulate the following:
    \begin{eqnarray}\label{maineq}
    \min_{\bf a_r,\bf b_r,\bf c_r,\bs \alpha_r}\Vert \mathcal{C} - \mathcal{L}\Vert_F^2 + \sigma \Vert \bs \alpha \Vert_{\ell_1}
    \end{eqnarray}
    where 
    $\mathcal{L} = \sum^{R}_{r=1} \bs \alpha_r \bf a_r \circ \bf b_r \circ \bf c_r$, $\Vert \cdot \Vert_F$ is the tensor Frobenius norm, and $\Vert \cdot \Vert_{\ell_1}$ is the vector one norm.
The optimization problem (\ref{maineq}) is neither convex nor differentiable \cite{KoBa09}. In addition, this formulation is amenable to tensor completion problems where one can find missing data. This is basis for predicting future Covid-19 infection cases in presence of previous pandemic data. 

Most of the algorithms suggested for solving (\ref{maineq}) are based on ALS (alternating least squares) \cite{KoBa09}. ALS is fast, easy to program and effective. However, ALS has some drawbacks. There is a need for more efficient and accurate methods. Thus, we propose a sampling method for ALS to leverage the efficiency of ALS and to allow bigger tensor data.

\subsection{Previous Work on Covid-19 Data Analysis}
Here, we mention a non-exhaustive collection of literature on Covid-19 analysis. The work varies from PDE models to machine learning methods for the prediction of Covid-19 hotspots, patterns, and outbreaks. A parabolic PDE-based predictive model with learned parameters from training data of preceding Covid-19 cases has been proposed to predict Covid-19 infections in Arizona\cite{pde_covid}. A tensor train global optimization is used to explore the parameter space to locate the starting points, and the Nelder-Mead simplex local optimization algorithm is used for the local optimization of the covid-19 modeling problem and eventually, Runge-Kutta method was applied to solve the PDE for one step forward prediction. There are several machine learning methods. In supervised machine learning, time series forecasting via Holt's Winter model was used to analyze global Covid-19 data and predict the sum of global Covid-19 cases in comparison to linear regression and support vector algorithms. Dictionary learning \cite{DBLP} through non-negative matrix factorization to identify the pattern and to predict future covid-19 outbreaks. In addition, machine learning has been applied to a susceptible-exposed-infected-removed (SEIR) model with SARS 2003 training data to predict covid-19 outbreaks \cite{Yang}. Deep learning models like attention-based encoder-decoder model \cite{mathew2020deep} have been implemented to forecast covid-19 outbreak. There needs to be more work on mathematical modeling to detect covid-19 hotspots; however, a tensor-based anomaly detection method for spatial-temporal data\cite{zhao2019rapid} has been proposed to determine hotspots based on an anomaly in pattern.
\subsection{Contributions.}
In this work, we focus on analyzing Covid-19 infection data \cite{CovidData,PopData} of New Jersey from the period of 04/01/2020 to 12/26/2021 (NJ-Covid-19). The state of New Jersey was initially chosen since we would like to study the spread of the disease in the most densely populated. The raw data collected by New York Times \cite{CovidData,PopData} was the daily basis cumulative data. Pre-processing techniques were applied to the raw spatio-temporal data and formatting the data in a tensor structure. Standard low-rank tensor decomposition models such as canonical polyadic (CP), Higher Order Orthogonal Iteration (HOOI), as well as tensor rank revealing methods like the tensor rank approximation method called low-rank approximation of tensor (LRAT) \cite{fat_nav} and LRAT with flexible Golub-Kahan \cite{nav_wang} are used to approximate the pattern and flow of covid infections. Novel stochastic canonical polyadic (STOCP) was also applied to the NJ-Covid-19 data. Some of the other contributions are the following:
\begin{itemize}
  
    \item 
    Converted each entry of  the Covid-19 tensor into the rate of increment in each week and applied the CP-ALS algorithm and plotted error with practical threshold using standard deviation and mean. The spike above the threshold line exactly replicates to the original spike of the covid-19 infections.
    \item The tensor completion optimization was formulated to predict future covid-19 infections as far as the following week and following quarter using LRAT with flexible Golub-Kahan.
    \item
    The proposed STOCP was implemented on the NJ-Covid-19 Data. To see its efficacy, we compare its numerical results with the standard CP decomposition model via the alternating least-squares method. We tested our spatiotemporal data tensor and other image data of different sizes.
\end{itemize}

\subsection{Outline of the paper}
The paper is organized as follows. In Section 2, we provide some tensor backgrounds, standard tensor decompositions, CP, Higher-Order Tucker, and well-known numerical techniques like Alternating Least-Squares and Higher-Order Orthogonal Iteration. Then, in Section 3, we include some discussions of the Covid-19 tensor as well as CP and HOT tensor model analysis. Section 4 describes a sampling method for ALS with some numerical results comparison with ALS.  Section 5 deals with the tensor sparse model which is implemented with LRAT with Golub-Kahan \cite{nav_wang}. We discuss its application to predicting Covid-19 infection cases and locating and identifying covid-19 hotspots. Finally, we provide concluding remarks and some future outlooks in Section 6.

\section{Preliminaries}
We denote a vector by a bold lower-case letter $\mathbf{a}$.
The bold upper-case letter $\mathbf{A}$ represents a matrix
and the symbol of a tensor is a calligraphic letter $\mathcal{A}$.
Throughout this paper, we focus on third-order tensors
$\mathcal{A}=(a_{ijk})\in\mathbb{R}^{I\times J\times K}$ of three indices $1\leq i\leq I,1\leq j\leq J$ and $1\leq k\leq K$,
but these can be easily extended to tensors of arbitrary order greater or equal to three.


The three kinds of matricization for third-order $\mathcal{A}$ are $\mathbf{A}_{(1)},\mathbf{A}_{(2)}$ and $\mathbf{A}_{(3)}$,  according to respectively arranging the column, row, and tube fibers to be columns of matrices, which are defined by fixing every index but one and denoted by $\mathbf{a}_{:jk}$, $\mathbf{a}_{i:k}$ and $\mathbf{a}_{ij:}$ respectively.
We also consider the vectorization for $\mathcal{A}$ to
obtain a row vector $\mathbf{a}$ such the elements of $\mathcal{A}$ are arranged
according to $k$ varying faster than $j$ and $j$ varying faster than $i$, i.e.,
$\mathbf{a}=(a_{111},\cdots,a_{11K},a_{121},\cdots,a_{12K},\cdots,a_{1J1},\cdots,a_{1JK},\cdots)$. Kronecker Product of two matrices $\mf A  \in \mb R^{m\times n}$ and $\mf B \in \mb R^{p\times q}$ is denoted as $\mf A\otimes \mf B \in \mb R^{mp \times nq}$ and obtained as the product of each element of $\mf A$ and the matrix$\mf B$. Khatri-Rao Product of two matrices $\mf A$ and $\mf B $ with same columns is  the column-wise Kronecker product:
    $
      \mathbf{A\odot B} = [\mf{a}_1\otimes\mf{b}_1\hspace{0.15cm}\hdots \hspace{0.15cm}\mf{a}_R\otimes \mf{b}_{R} ]. 
   $
The outer product of a rank-one third order tensor is denoted as $\mathbf{a}\circ\mathbf{b}\circ\mathbf{c}\in\mathbb{R}^{I\times J\times K}$ of three nonzero vectors $\mathbf{a}, \mathbf{b}$ and $\mathbf{c}$ is a rank-one tensor with elements $a_ib_jc_k$ for all the indices, i.e. the matricizations of  $\mathbf{a}\circ\mathbf{b}\circ\mathbf{c}\in\mathbb{R}^{I\times J\times K}$ are rank-one matrices.
A canonical polyadic (CP) decomposition of $\mathcal{A}\in\mathbb{R}^{I\times J\times K}$ expresses $\mathcal{A}$ as a sum of rank-one outer products:
\begin{equation}\label{cpd}
\mathcal{A}=\sum_{r=1}^{R} \mathbf{a}_r\circ\mathbf{b}_r\circ\mathbf{c}_r
\end{equation}
where $\mathbf{a}_r\in\mathbb{R}^I,\mathbf{b}_r\in\mathbb{R}^J,\mathbf{c}_r\in\mathbb{R}^K$ for $1\leq r\leq R$ and and $\circ$ is the outer product.
The outer product $\mathbf{a}_r \circ \mathbf{b}_r \circ \mathbf{c}_r$ is a rank-one component
and the integer $R$ is the number of rank-one components in tensor $\mathcal{A}$.
The minimal number $R$ such that the decomposition (\ref{cpd}) holds is the rank of tensor $\mathcal{A}$, which is denoted by $\mbox{rank}(\mathcal{A})$.
For any tensor $\mathcal{A}\in\mathbb{R}^{I\times J\times K}$,
$\mbox{rank}(\mathcal{A})$ has an upper bound $\min\{IJ,JK,IK\}$ \cite{kruskal}. In fact, tensor rank is NP-hard over $\mathbb{R}$ and $\mathbb{C}$ \cite{10.1145/2512329}.
Another standard tensor decomposition is higher-order Tucker (HOT) decomposition. HOT is a generalization of matrix SVD where $m \times n$ matrix $\mathbf{M}$ has a factorizaion $\mathbf{U} \mathbf{\Sigma} \mathbf{V}^T$, where $\mathbf{U}$ and $\mathbf{V}$ are orthogonal matrices and $\Sigma$ is a diagonal matrix with $\sigma_{ij}=0$ if $i \neq 0$ and otherwise $\sigma_{ii}>=0$. Its generalization in third-order tensors is
\begin{eqnarray*}
\mathcal{T}=\mathcal{G}\times_{1}\mathbf{U}^{(1)}\times_{2}\mathbf{U}^{(2)}\times_{3}\mathbf{U}^{(3)}
\end{eqnarray*}
where $\mathcal{T}\in \mathbb{R}^{I_1 \times I_2 \times I_3}$ is the given tensor, $\mathcal{G}\in \mathbb{R}^{R_1 \times R_2 \times R_3}$ is the core tensor and 
  $\mathbf{U^{(i)}} \in \mathbb{R}^{I_i \times R_i}$ for $i=1,2,3$ is an orthogonal matrix. The Tucker contracted product $\times_{1}$ is defined as $\mathcal{G}\times_{1}\mathbf{U}^{(1)}= \sum_{r_1}\mathcal{G}_{r_1 r_2 r_3} \mathbf{U}^{(1)}_{i_1 r_1} \in \mathbb{R}^{I_1 \times R_2 \times R_3}$.

\subsection{Standard Least-Squares Optimization for Tensor Decomposition}
Tensor decompositions like CP and HOT  \cite{BRO1997149,KoBa09} are considered to be generalizations of the singular value decomposition (SVD) and principal component analysis (PCA) of a matrix. To achieve CP from a given third order $\mathcal{X}$, an optimization problem is solved to find the finite sum of rank one tensor of $3^{rd}$ order approximating $\mathcal{X}$:
\begin{align}
    \min_{\mathbf{A,B,C}} \mid\mid \mathcal{X} - \mathcal{\hat{X}} \mid \mid_F
\end{align}    
\begin{center}
    with
\end{center}
\begin{align}
    \mathcal{\hat{X}} = \sum^{R}_{r=1}\mathbf{a_r}\circ\mathbf{ b_r }\circ \mathbf{c_r }
\end{align}
where $\mathbf{A,B,C}$ are factor matrices with $\mathbf{a_r,b_r,c_r}$ are their respective column vectors. This nonlinear optimization can be divided into subproblems of linear least squares with respect to the factor matrices. This is called the Alternating Least Squares (ALS).

\begin{figure}[H]
    \centering
    \includegraphics[scale = 1.2]{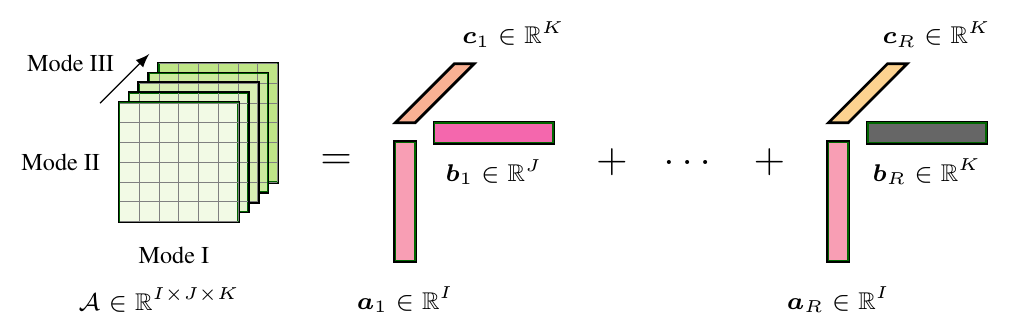}
    \caption{\textbf{CP-Decomposition Architecture}}
    \label{fig:cp_arch}
\end{figure}

ALS is an iterative method for finding the CP decomposition of a given tensor.
The nonlinear optimization problem is:
\begin{align}
\min_{\mathbf{A,B,C}}\vert\vert \mathcal{X} - \sum_{r= 1}^{R}\mathbf{a_r\circ b_r \circ \ c_r}\vert\vert_{F} ^2
\end{align}
where  $\mathbf{A}$, $\mathbf{B}$ and $\mathbf{C}$ are factor matrices containing the columns $\mathbf{a_r} \in \mathcal{R}^I$, $\mathbf{b_r} \in \mathcal{R}^J$ and $\mathbf{c_r} \in \mathcal{R}^K$ respectively.
The problem can be reduced to linear least squares problems at each iteration with an initial guess $\mathbf{A^0,B^0,C^0}$, the sequences $\mathbf{A^k,B^k,C^k}$ are generated by solving each sub-problems\cite{3831,lieven}.  Given the initial guess $\mathbf{A}^0$,$\mathbf{B}^0$,$\mathbf{C}^0$, the factor matrices are updated by the following:
update $\mathbf{A}$ via
\begin{align}
    \mathbf{A}^{k+1} = \argmin_{\mathbf{A} \in \mathcal{R}^{I\times R}}{\frac{1}{2}\mid\mid \mathbf{X_{(1)} - \mathbf{A(C^k\odot B^k)^T}}\mid\mid}^2_F
\end{align}
update $\mathbf{B}$ via
\begin{align}
    \mathbf{B}^{k+1} = \argmin_{\mathbf{B} \in \mathcal{R}^{J\times R}}{\frac{1}{2}\mid\mid \mathbf{X_{(2)} - \mathbf{B(C^k\odot A^k)^T}}\mid\mid}^2_F
\end{align}
and update $\mathbf{C}$ via
\begin{align}
    \mathbf{C}^{k+1} = \argmin_{\mathbf{C} \in \mathcal{R}^{K\times R}}{\frac{1}{2}\mid\mid \mathbf{X_{(3)} - \mathbf{C(B^k\odot A^k)^T}}\mid\mid}^2_F
\end{align}
These updating schemes are repeated until convergence; See figure \ref{fig:cp_arch}, Algorithm \ref{alg:cap}. The local convergence and uniqueness of the ALS technique were discussed in this work \cite{local_conv}.


 \begin{algorithm}[H]
\caption{\cite{Nicholas}CP-ALS for 3-Way Tensor Decomposition}\label{alg:cap}
 \begin{algorithmic}
\Require Tensor $\mathcal{A} \in \mathbb{R}^{I \times J \times K}$rank R,Maximum Iterations N.
\Ensure CP Decomposition $\mathbf{\lambda}\in \mathrm{R}^{R \times 1}$, $\mathbf{A} \in \mathrm{R}^{I \times R}$,$\mathbf{B} \in \mathrm{R}^{J \times R}$, $\mathbf{C} \in \mathrm{R}^{K \times R}$.
\\
1: Initialize $\mathbf{A,B,C};$\\
2: \textbf{for} i = 1 ... N \textbf{DO}\\
3:\hspace{0.7cm} $\mathbf{A} \gets \mathbf{X}_{(1)}(\mathbf{C\odot B})(\mathbf{C^{T}C*B^{T}B})^{\dagger}$\\
4:\hspace{0.7cm} Normalize the columns of $\mathbf{A}$(storing norms in vector $\lambda.)$\\
5:\hspace{0.7cm} $\mathbf{B} \gets \mathbf{X}_{(2)}(\mathbf{C\odot A})(\mathbf{C^{T}C*A^{T}A})^{\dagger}$\\
6:\hspace{0.7cm} Normalize the columns of $\mathbf{A}$(storing norms in vector $\lambda.)$\\
7:\hspace{0.7cm} $\mathbf{C} \gets \mathbf{X}_{(3)}(\mathbf{B\odot A})(\mathbf{B^{T}B*A^{T}B})^{\dagger}$\\
8:\hspace{0.7cm} Normalize the columns of $\mathbf{A}$(storing norms in vector $\lambda.)$\\
9:\hspace{0.7cm}\textbf{If} convergence is met \textbf{then},\\
10:\hspace{1.2cm}\textbf{break for loop}\\
11:\hspace{0.7cm}\textbf{end if}\\
12: \textbf{end for}\\
13: return $\mathbf{\lambda,A,B,C}$;
\end{algorithmic}

\end{algorithm}

 \begin{algorithm}[H]
\caption{The CP-ALS Algorithm for General Case}
 \begin{algorithmic}
\Require Initialize with factor matrics $ \mathbf{A^0} \in 
\mathbb{R}^{I\times R},\mf B^0\in \mathbb{R}^{J\times R},\mf C^0 \in \mathbb{R}^{K\times R}$
\Ensure Factor matrices $\mathbf{A^{(n)}}\in \mb R^{I_n \times R}$ for n = $1,2,\hdots,N.$
\\
1: \textbf{for} i = k ... N \textbf{DO}\\
2:\begin{align*}
  (\mf A^{(n)})^{k+1}& = \argmin_{\mf A^{(n)}}\frac{1}{2}\vert\vert \mf X_{(n)} - \mathbf{A((A^{(N)})^k\odot \hdots\odot(A^{(n + 1)})^k\circ (\mf A^{(n-1)}}^{k+1}\\
  & \hspace{5cm} \odot \hdots \odot(\mf A^{(1)})^{k+1})^T\vert\vert_F^2 
\end{align*}
3: \textbf{end for}\\
4: return $\mathbf{A^{(n)}}\in \mb R^{I_n \times R}$;
\end{algorithmic}
\end{algorithm}

\subsection{Higher Order Orthogonal Iteration (HOOI)}
\begin{figure}[H]
      \centering
      \includegraphics[scale = 0.3]{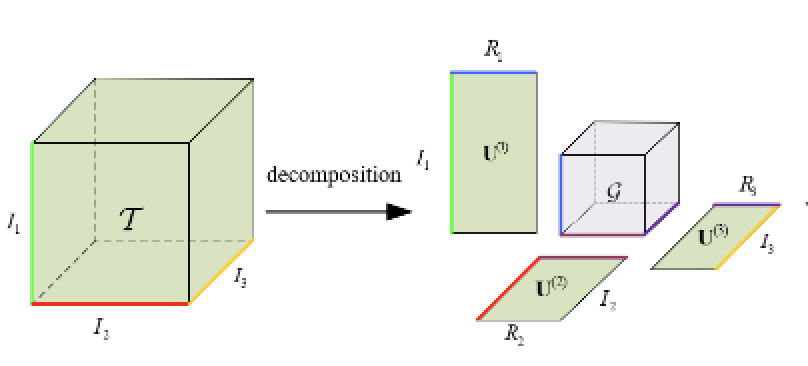}
      \caption{HOOI Architecture of 3-Way Tensor\cite{e23101349}}
      \label{fig:my_label2}
  \end{figure}

HOOI is an iterative algorithm that computes low rank decomposition of a given tensor \cite{HOOI}. To achieve the reconstruction of $\mathcal{T}$ in the HOOI format, i.e. $\mathcal{T} = \mathcal{G} \times_{1} \mathbf{U}^{(1)}\times_{2}\mathbf{U}^{(2)}... \times_{N}\mathbf{U}^{(N)}$, the problem formulation is
  \begin{align*}
      \min_{\mathbf{U^{(1)},U^{(2)},...,U^{(N)}}}{\mid\mid \mathcal{T}-\mathcal{G}\times_{1}\mathbf{U}^{(1)}\times_{2}\mathbf{U}^{(2)}... \times_{N}\mathbf{U}^{(N)}\mid\mid_{F}^{2}}
  \end{align*}
  where $\mathcal{T}\in \mathbb{R}^{I_1 \times ... \times I_N}$ and 
  $\mathbf{U^{(i)}} \in \mathbb{R}^{I_i \times R_i}$ is a low rank matrix for $i=1,\ldots, N$ with $I_i \geq R_i$. Equivalently, the optimization can be reformulated as a maximization of $\mid\mid \mathcal{G}\mid\mid_{F}^{2}$ with $\mathcal{G}$ is given by
  \begin{align}
       \mathcal{G} = \mathcal{A}\times_{1}\mathbf{U}^{(1)^T}\times_{2}\mathbf{U}^{(2)^T}... \times_{N}\mathbf{U}^{(N)^T}
  \end{align}   
where the core tensor $\mathcal{G} \in \mathbb{R}^{R_1 \times \ldots R_N}$ tensor for $i=1,\ldots,N$. See Figures \ref{fig:my_label2}.
  \begin{algorithm}[H]
\caption{\cite{fusion}Higher-Order Orthogonal Iteration(HOOI)}
  \begin{algorithmic}
\Require :Tensor $\mathcal{X} \in \mathbb{R}^{I_{1}\times ... \times I_{N}}$ and ranks $R_{1},...,R_{N}$
\Ensure:Tucker Factors:$\mathbf{U_{1} \in \mathbb{R}^{I_{1}\times R_{1}},U_{2}\in \mathbb{R}^{I_{2}\times R_{2}},...,U_{N} \in \mathbb{R}^{I_{N}\times R_{N}}}$\\and Core tensor $\mathcal{G} \in \mathbb{R}^{R_{1}\times R_{2}\times ... \times R_{N}}$\\
1: Initialize $\mathbf{U_{1},...,U_{N}};$(random or given by HSVD)\\
2: \textbf{while} convergence criterion not met \textbf{do}\\
3: \hspace{0.6cm}\textbf{for} n = 1,...,N \textbf{do}\\
4: \hspace{0.9cm} $\mathcal{W} \gets \mathcal{X}\times_{N}\mathbf{U_{N}^{T} ... \times_{n+1} U_{n+1}^{T}\times_{n-1}U_{n-1}^T ... \times_{1}U_{1}^{T}}$\\
5: \hspace{0.9cm} $\mathbf{[U \Sigma V]} \gets SVD(W_{(n)})$\\
6: \hspace{0.9cm} $\mathbf{U_{n} \gets U(:,1:R_n)}$\\
7: \hspace{0.6cm}\textbf{end for}\\
8: \textbf{end while}\\
9: $\mathcal{G} \gets \mathcal{X} \times_{N}\mathbf{U_{N}^{T}\times_{N-1}U_{N-1}^{T} ... \times_{1}U_{1}^{T}}$
\end{algorithmic}
\end{algorithm}

\section{Covid-19 Tensor Analysis Using ALS and HOOI}
\subsection{Covid-19 Data Preparation}
We focus on the daily Covid-19 infection data of New Jersey State from the New York Times's Covid-19 Data Depository from the period of 04/01/2020 to 12/26/2021 \cite{CovidData,PopData}. The state of New Jersey was initially chosen since we would like to investigate the spread of the disease in the most densely populated and affected state. The raw data collected by New York Times was a daily basis cumulative data. We restructured the data table into weekly total Covid-19 infections in each county. We stacked 21 matrices representing the counties of New Jersey  of size 13 weeks $\times$ 7 quarters. Thus, we construct a tensor data, $\mathcal{C} \in \mathbb{R}^{13 \times 7 \times 21}$. Each element of the tensor represents the total infection in a week of a particular quarter in a county. See Figure \ref{fig:covid-tensor} below.
\begin{figure}[H]
\centering
\begin{tikzpicture}

\coordinate (O) at (0, 0, 0);
\coordinate (A) at (0, \Width, 0);
\coordinate (B) at (0, \Width, \Height);
\coordinate (C) at (0, 0, \Height);
\coordinate (D) at (\Depth, 0, 0);
\coordinate (E) at (\Depth, \Width, 0);
\coordinate (F) at (\Depth, \Width, \Height);
\coordinate (G) at (\Depth, 0, \Height);

\draw[green!60!black, fill = green!5] (O) -- (C) -- (G) -- (D) -- cycle; 
\draw[green!60!black, fill = green!5] (O) -- (A) -- (E) -- (D) -- cycle; 
\draw[green!60!black, fill = green!5] (O) -- (A) -- (B) -- (C) -- cycle; 
\draw[green!60!black, fill = green!5, opacity = 0.2] (D) -- (E) -- (F) -- (G) -- cycle; 
\draw[green!60!black, fill = green!5, opacity = 0.4] (C) -- (B) -- (F) -- (G) -- cycle; 
\draw[green!60!black, fill = green!5, opacity = 0.4] (A) -- (B) -- (F) -- (E) -- cycle; 

\draw (0.2, -1.2, 0) node[rotate = 0,scale =0.6 ] {\scriptsize{\color{purple}$7$~Days}};
\draw (0.2, -0.95, 0) node[rotate = 0] {{\color{purple!65}$\underbrace{\hspace{2cm}}$}};
\draw (-0.4, 1, 2) node[rotate = 90,scale = 0.5] {\scriptsize{\color{purple}$13$~Weeks}};
\draw (-0.15, 1, 2) node[rotate = 270] {{\color{purple!65}$\underbrace{\hspace{2cm}}$}};
\draw (2.5, -0.2, 1.4) node[rotate = 45,scale = 0.5] {\scriptsize{\color{purple}$21$~Counties}};
\draw (2.2, 0, 1.2) node[rotate = 45] {{\color{purple!65}$\underbrace{\hspace{1.1cm}}$}};
\draw (1.3, 1.4, 2.) node[rotate = 0,scale =2] {\scriptsize{\color{black}$\mathcal{C}$}};
\end{tikzpicture}
\caption{Covid-19 Tensor: New Jersey's number of weekly Covid-19 cases every quarter per county.}
\label{fig:covid-tensor}
\end{figure}
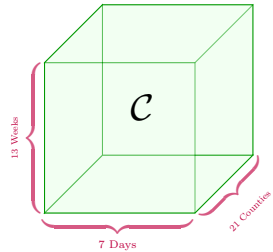
To increase the accuracy and efficacy of the algorithms, we normalized the data tensor converting it into relative cases tensor with respect to the population of the respective county. We constructed the population tensor of New Jersey collecting population data from US Census 2020 \cite{PopData} and divided the tensor-data tensor by the population data tensor. Each element of our new normalized tensor is the following:
\begin{align*}
        \mf C(i,j,k)& = \frac{\text{Total infections in } j^{th}\text{ week of } i^{th} \text{ quarter in the } k^{th} \text{ county}}{\text{Respective Mid-year Population}}\\
\end{align*}
Specifically, we have
\begin{align*}
        \mf C(1,1,1) & = \frac{\text{Total covid infections of Atlantic County in 1st week of April}}{\text{mid year population of Atlantic county in 2020}}\\
        &= \frac{\mf A(1,1,1)}{\mf P(1,1,1)}= \frac{67}{274534}= 2.44e-04
    \end{align*}

The tensor $\mathcal{C}$ is rescaled by dividing by the population of respective counties of mid-year 2020 \cite{PopData}. The new tensor data $\bar{\mathcal{C}}$ has greatly improved efficiency in the numerical implementations. See sections \ref{cp-hooi},\ref{pattern-factors},\ref{sto-als}, and \ref{tc-lrat} for some experimental results. 

\subsection{Numerical Experiments}
\subsubsection{CP vs HOOI}
\label{cp-hooi}
\begin{figure}[H]
     \centering
     \includegraphics[scale = 0.29]{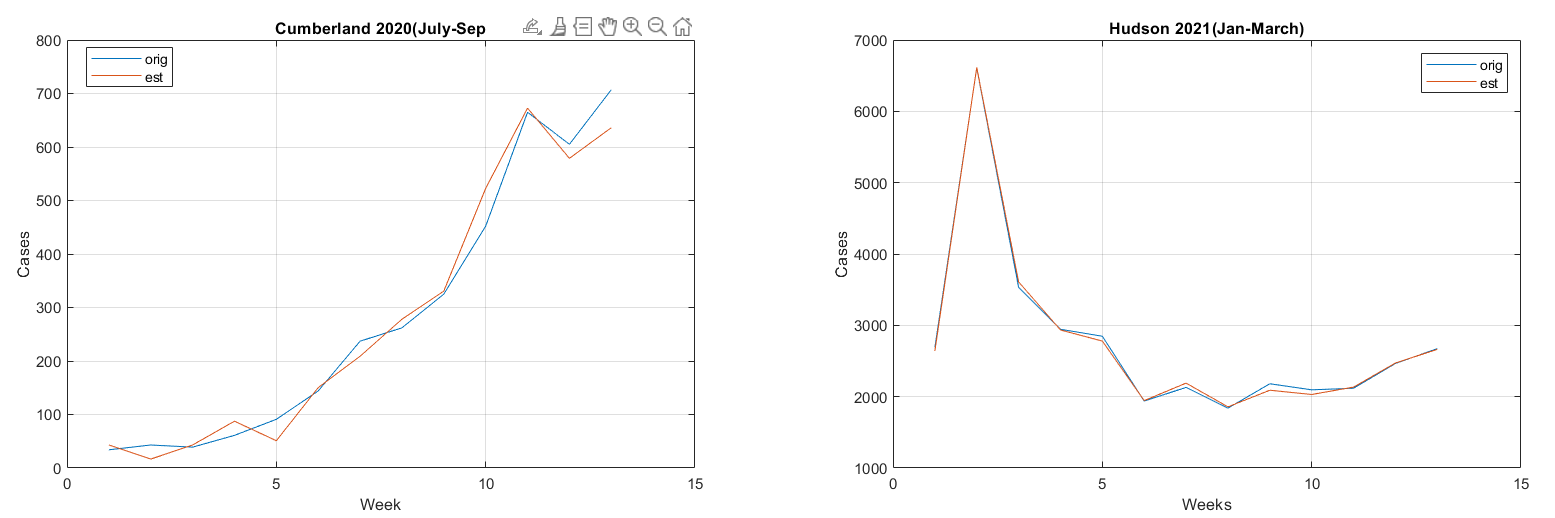}
     \caption{Original and Approximated using CP}
     \label{fig:my_label-1}
 \end{figure}

\begin{figure}[H]
     \centering
     \includegraphics[scale = 0.3]{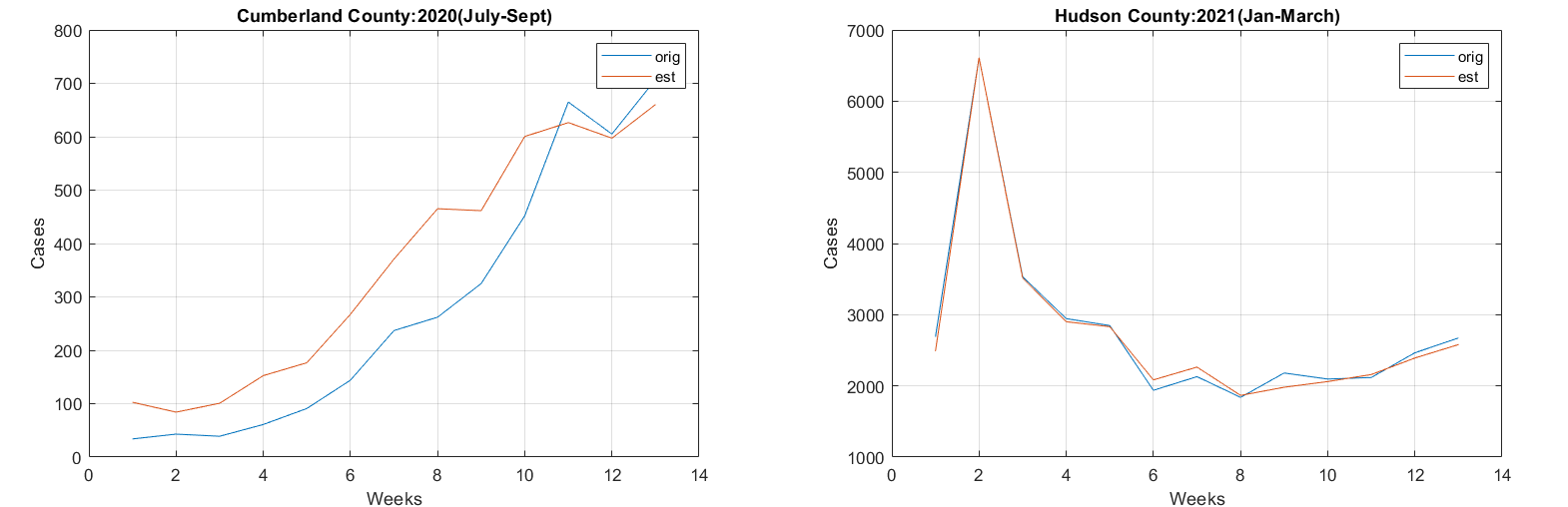}
     \caption{Original and Approximated using HOOI}
     \label{fig:my_label-2}
 \end{figure}

\begin{figure}[H]
    \centering
    \includegraphics[scale= 0.3]{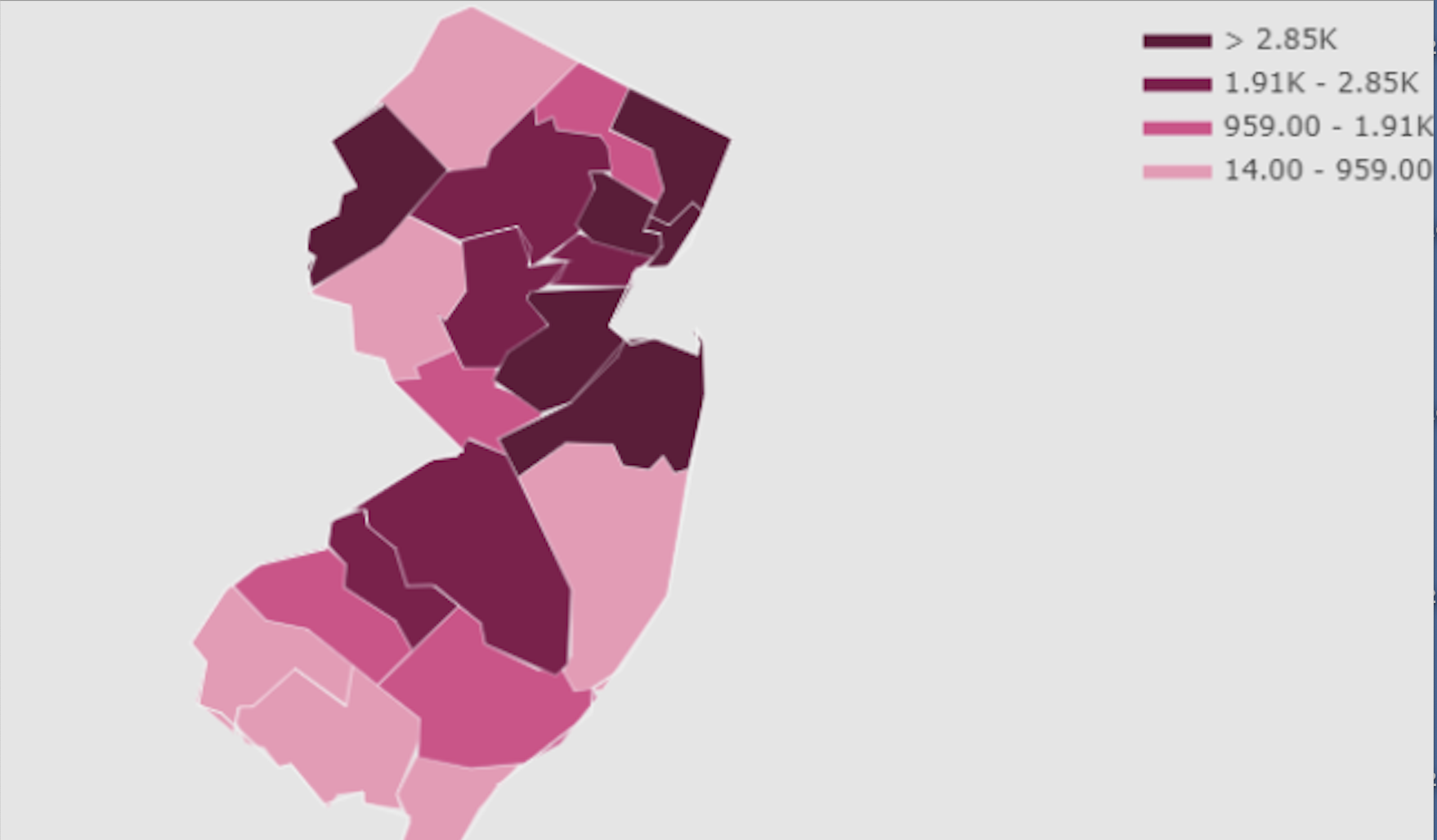}
    \caption{Original Covid-19 Cases Geo plot: 3rd
Week of Jan 2021}
    \label{fig:my_label-3}
\end{figure}

\begin{figure}[H]
    \centering
    \includegraphics[scale= 0.3]{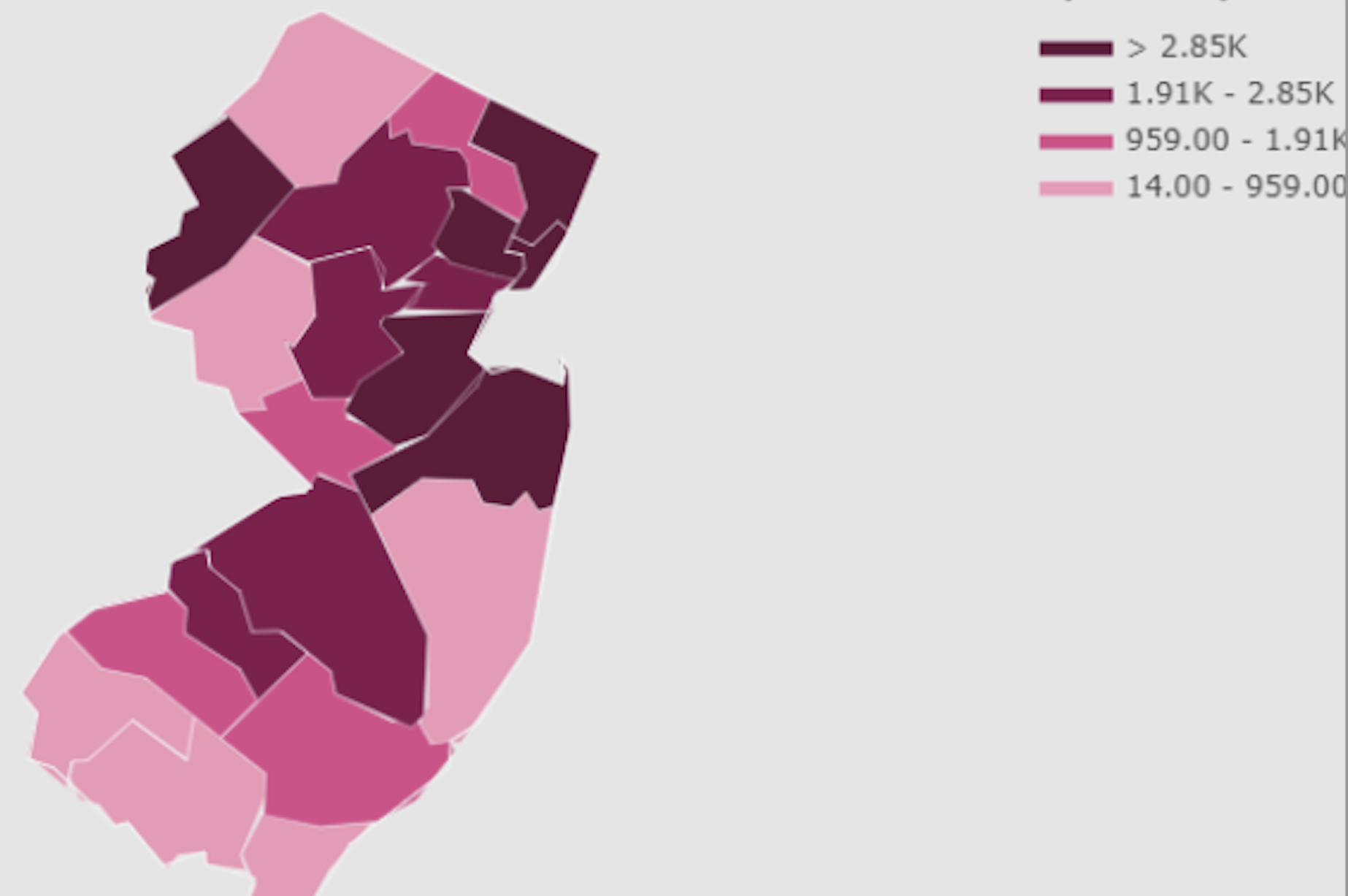}
    \caption{Approximated Lower Rank CP with R
= 20, 3rd Week Jan 2021}
    \label{fig:my_label-4}
\end{figure}

One main advantage of tensor decomposition, namely, ALS and HOOI, is that it provides analytic tools for higher-order data in several modes. We implemented the algorithms, CP and HOOI, on our tensor $\bar{\mathcal{C}}$ of size $13 \times 7\times 6$. First, we ran the ALS algorithm to construct the CP decomposition into three-factor matrices. Then we analyze through the visualization of the three-factor matrices from the estimated constructed tensor. We are able to estimate the same evolutionary patterns of Covid-19 cases as the original data tensor; see Figures \ref{fig:my_label-1} and \ref{fig:my_label-2}. Collectively, we plot the cases in Figures \ref{fig:my_label-3} and \ref{fig:my_label-4} based on the intensity levels. HOOI is faster in time but it gives lower accuracy results than CP.
 \begin{center}
    \textbf{Geo-Plot of CP approximation}
\end{center}
\subsubsection{Extracting Patterns in Covid-19 via Factor Matrices}
\label{pattern-factors}
\begin{figure}[H]
    \centering
    \includegraphics[scale = 0.4]{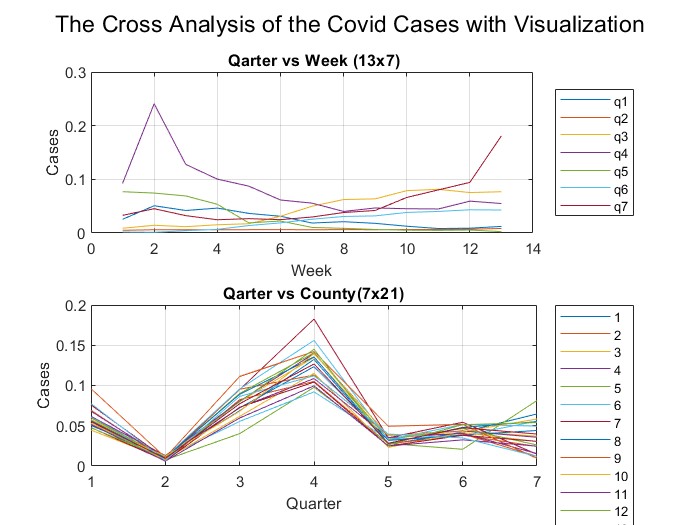}
    \caption{Cross-Pattern on Sep 27-Dec-26 2020  }
    \label{fig:my_label7}
\end{figure}

\begin{figure}[H]
    \centering
    \includegraphics[scale = 0.65]{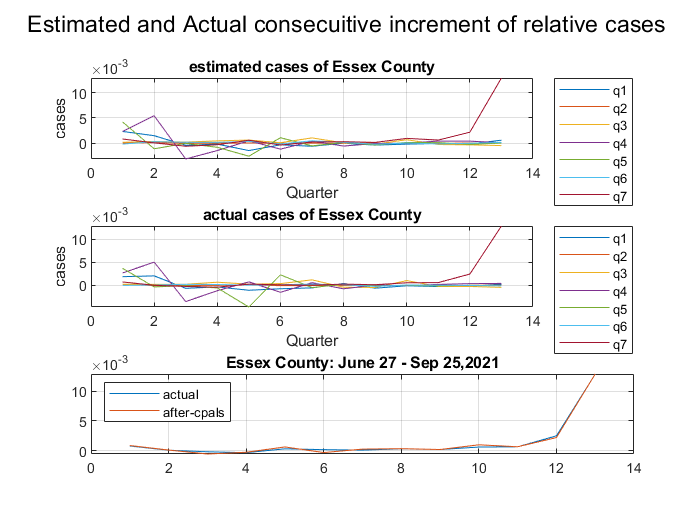}
    \caption{This is the Covid-19 weekly consecutive increment (relative cases) of Essex county.7th quarter of June 27 - Sept 25, 2021. }
    \label{fig:my_label8}
\end{figure}

We further explore from our output factor matrices of the CP decomposition. We multiply different factor matrices to observe the county-wise, week-wise, and quarter-wise pattern of the Covid-19 cases. The first visualization of Figure \ref{fig:my_label7} shows which quarter has the highest cases increment in 2nd week. Then the second visualization of the quarter, Sep 27-Dec-26 2020, indicates that the fourth week has the highest number of cases. Our algorithm estimates that the week, of Oct 4-10,2020, is the week of the highest increment. Furthermore, we found that the same identical patterns on the original Covid-19 data tensor are consistent with our findings. Similarly, Figure \ref{fig:my_label8} shows the identical pattern of our estimated and original Covid-19 infections in Essex county from June 27 - Sept 25, 2021. 



\section{Sampling Method for Alternating Least-Squares (SMALS)}
Given  a tensor $\mathcal{X} \in \mathbb{R}^{I\times J \times K}$ and a fixed positive integer $R$, the ALS algorithm solves three independent least squares problems alternatively. The least squares problems can be solved via various methods such as QR factorization, Cholesky factorization and etc. However, it requires of an $\mathcal{O}(R^3)$ floating point operations per second. To reduce this we let $S \subset \{ 1, \hdots R\}$ be the set of sample indices and $A_s,  B_s$ and $C_s$ represent the sub-matrices obtained by choosing the columns of $A, B$, and $C$ according to the index set $S$ respectively. The partial derivatives of the objective function $f$ with respect to the blocks $A_s, B_s$ and $C_s$ are 
\begin{equation}
    \frac{\partial f}{\partial A_S} = -X_{(1)}(C_S \odot B_S)+A(C \odot B)^T(C_S \odot B_S),
\end{equation}
\begin{equation}
    \frac{\partial f}{\partial B_S} = -X_{(2)}(C_S \odot A_S)+B(C \odot A)^T(C_S \odot A_S),
\end{equation}
and 
\begin{equation}
    \frac{\partial f}{\partial C_S} = -X_{(3)}(B_S \odot A_S)+C(B \odot A)^T(B_S \odot A_S).
\end{equation}
the stationary points of the above equations can be obtained by setting each gradient equal to zero. For instance  $\nabla_{A_S} f =0$ implies the following normal equation
\begin{align}
    A_S\left( (C \odot B)^T (C_S \odot B_S) \right)^S &= -A_{S^C} \left( (C \odot B)^T (C_S \odot B_S) \right)^{S^C}\\ &+ X_{(1)}(C_S \odot B_s) \nonumber
\end{align}
where $A^S$ represents the sub matrix of $A$ obtained by sampling the rows of $A$ corresponding to the sampling set $S$. Similar results can be derived by solving the equations $\nabla_{B_S} f=0$ and $\nabla_{C_S}f=0$. This reduces the latter computational complexity to $\mathcal{O}(\max \{\vert S\vert^3\})$. When $I, J\, \text{and}\, K$ are relatively large, the reduction could be significant.  The sampling set $S$ is chosen based on the performance of each block variable in each iteration. For instance, if the update for the block $A_i$ yields more decrease in the objective function compared to the block $A_j$, the index $i$ is replaced by $j$ in the next iteration. The differentiation of ALS can be found here \cite{7891546}. The derivation of SMALS and further analysis can be found \cite{sampling_cp}.


\subsection{Numerical Results of SMALS vs ALS}
\label{sto-als}
We generate some numerical results to show the efficacy of SMALS and to compare with standard ALS on Covid-19 tensor data and random color image data.
\begin{figure}[H]
    \centering
    \includegraphics[scale = 0.17]{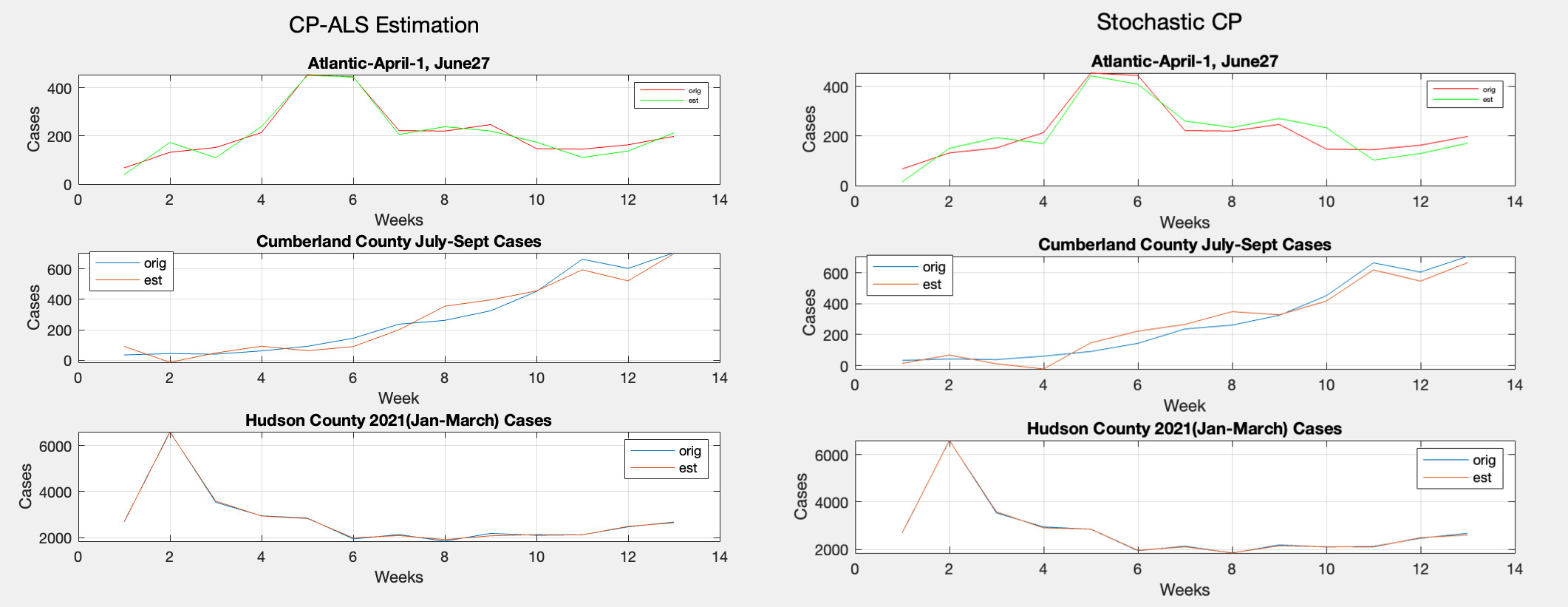}
    \caption{SMALS vs ALS on Real Covid-19 Tensor}
    \label{fig:my_label3b}
\end{figure}
\begin{figure}[H]
    \centering
    \includegraphics[scale = 0.17]{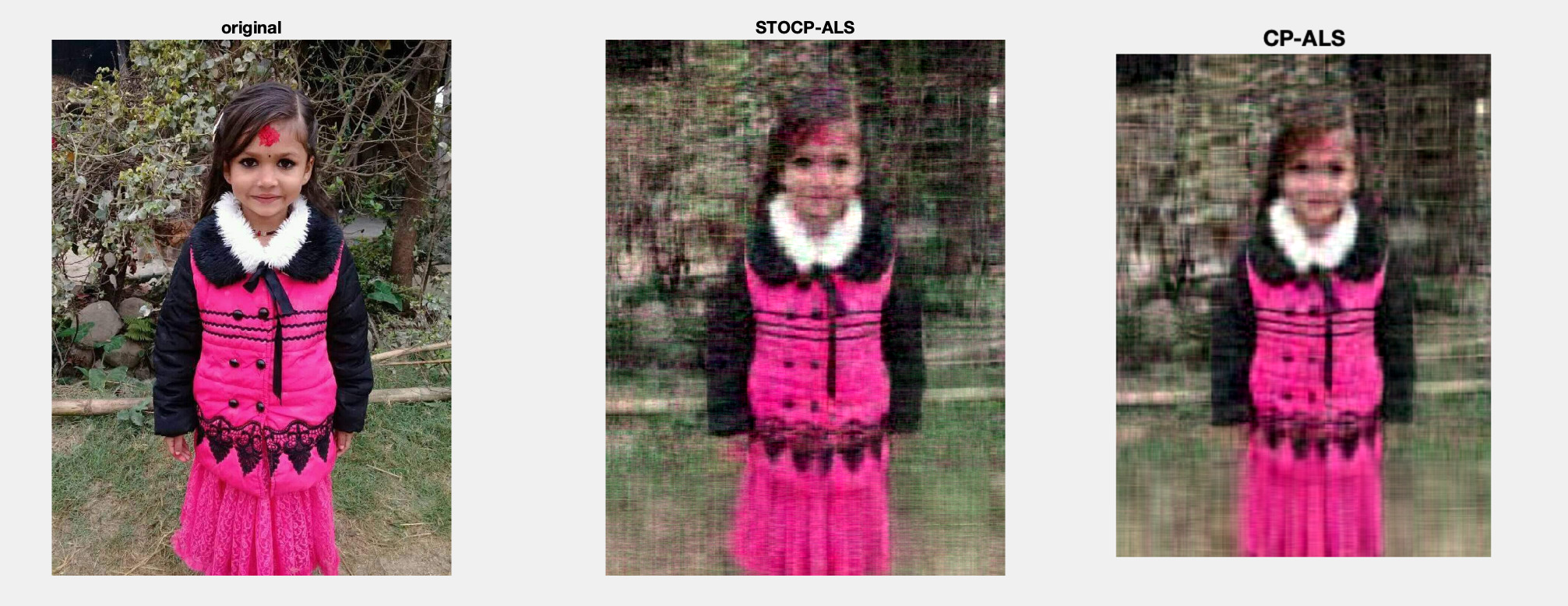}
    \caption{SMALS vs ALS on Random Color Images}
    \label{fig:my_label3c}
\end{figure}
\begin{figure}[H]
    \centering
    \includegraphics[scale = 0.17]{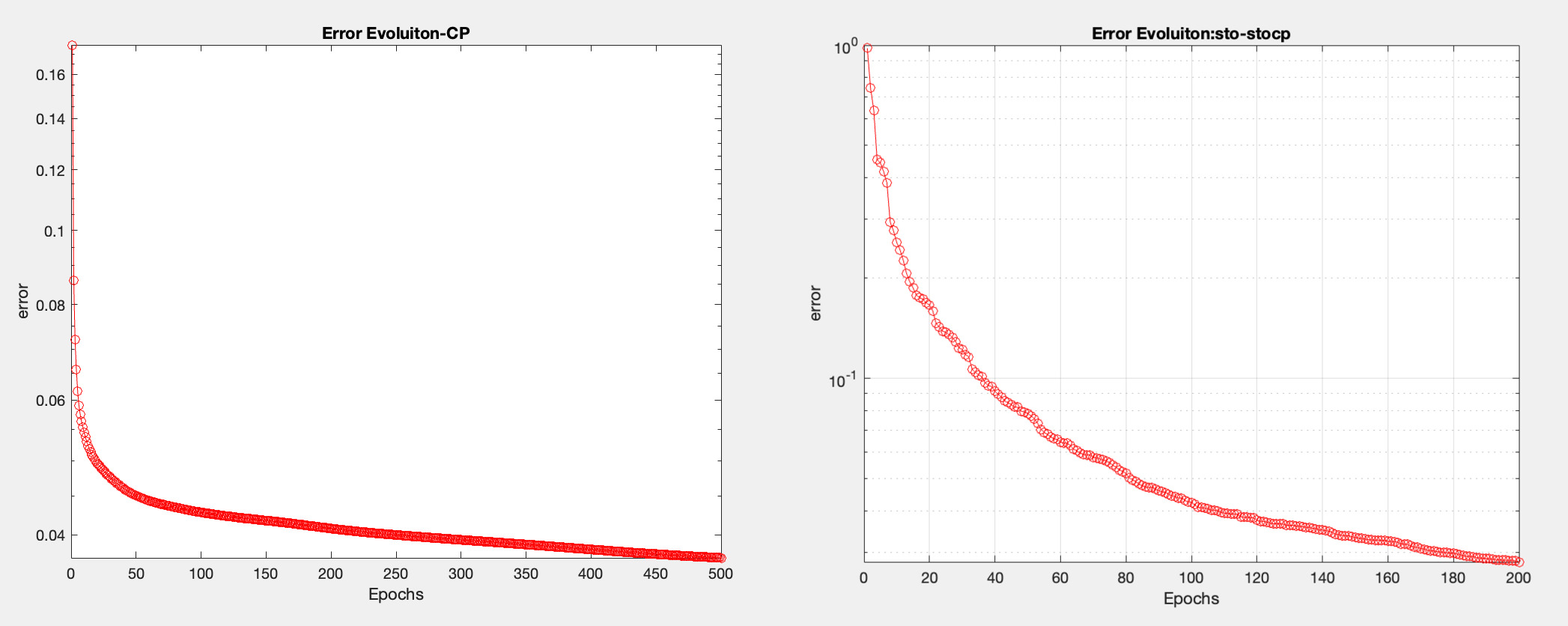}
    \caption{Error Evolution of ALS vs SMALS on random images in fig -8}
    \label{fig:my_label3d}
\end{figure}
\begin{center}
    \textbf{Efficacy Comparison Between ALS and SMALS}
\end{center}
\begin{figure}[H]
    \centering
    \includegraphics[scale = 0.3]{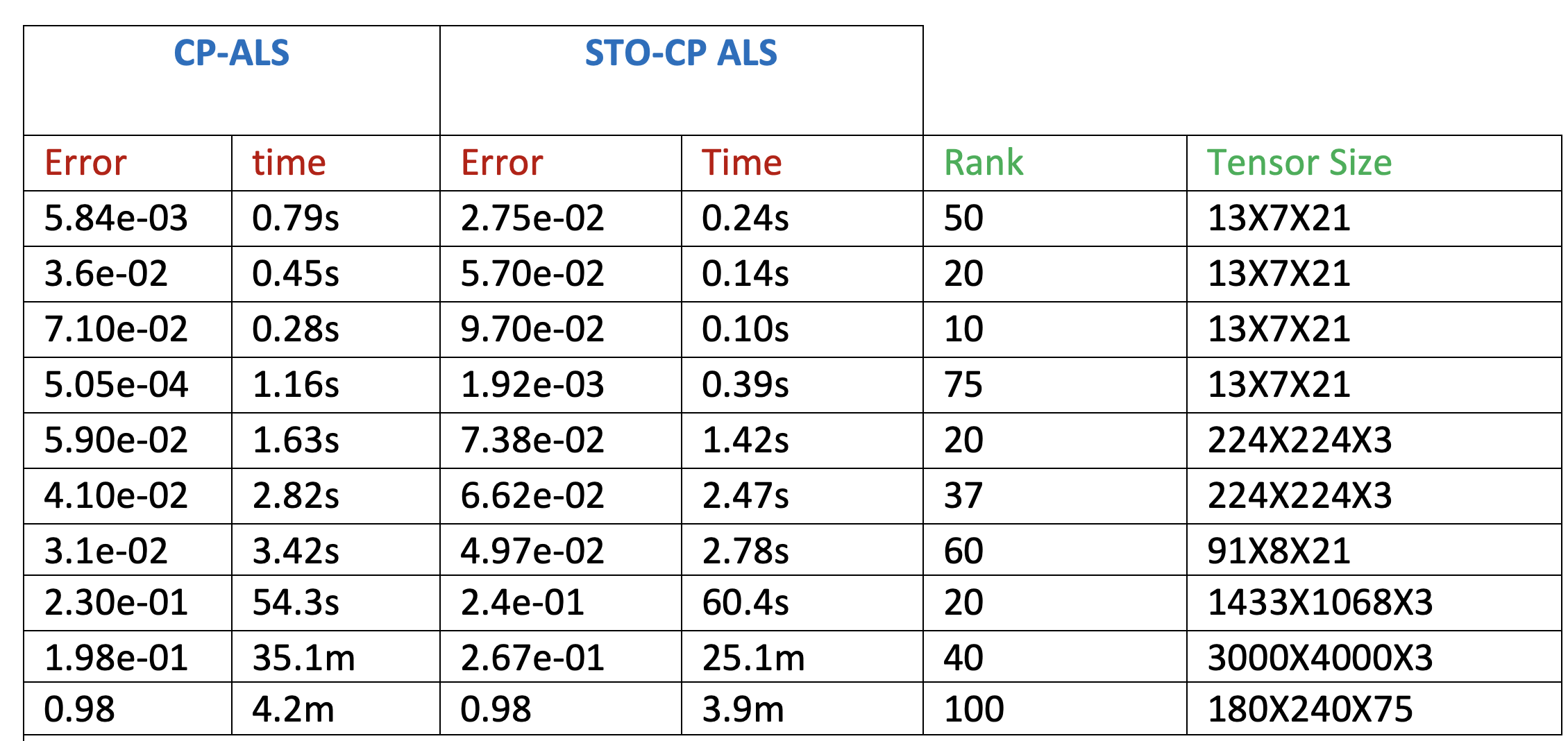}
    \caption{ALS vs SMALS Comparison}
    \label{fig:my_label3e}
\end{figure}
We implement both CP-ALS and CP-SMALS algorithms on different sizes of varying tensors. We run the codes twenty times for each tensor case and average their results. Our result indicates SMALS is more efficient in terms of time cost, see figure \ref{fig:my_label3e}

\section{Tensor Sparse Optimization}
\begin{figure}[H]
    \centering
    \includegraphics[scale = 0.6]{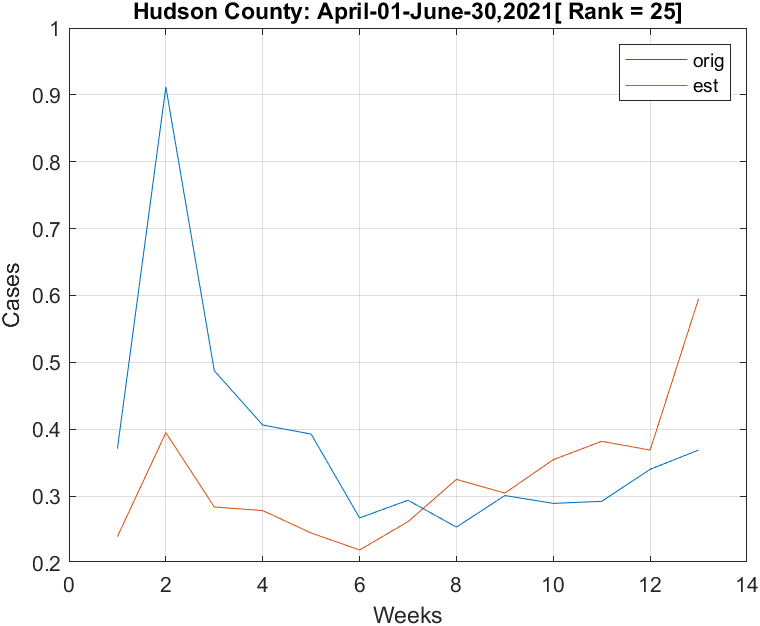}
    \caption{LRAT + fGK Estimation}
    \label{fig:my_label4}
\end{figure}

\begin{figure}[H]
    \centering
    \includegraphics[scale = 0.6]{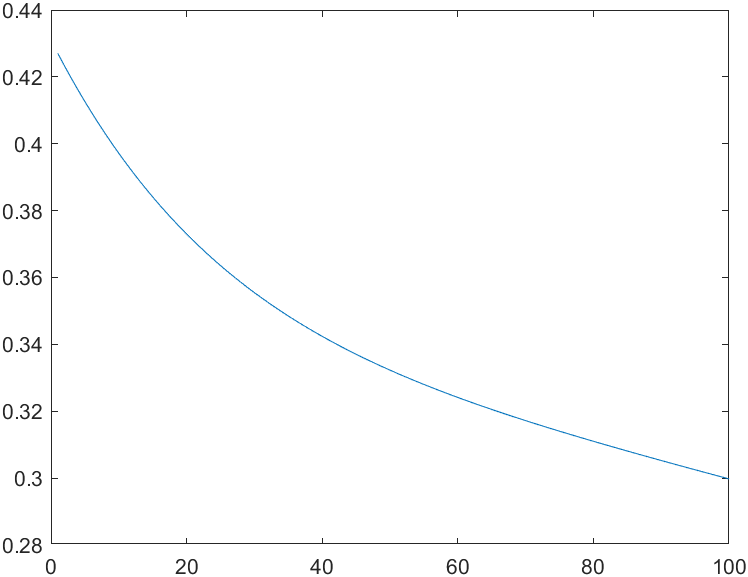}
    \caption{LRAT + fGK Convergence Plot}
    \label{fig:my_label5}
\end{figure}

In \cite{nav_wang}, an iterative method based on proximal algorithms called low-rank approximation of tensors (LRAT) was proposed to solve the minimum rank optimization:
\begin{eqnarray} \label{tensor_sparse_model}
    \min_{a_r,b_r,c_r,\alpha_r}\Vert \mathcal{C} - \mathcal{L}\Vert_F^2 + \sigma \Vert \bm\alpha \Vert_{\ell_1}
\end{eqnarray}
    where 
    $\mathcal{L} = \sum^{R}_{r=1}\alpha_r \mathbf{a_r} \circ \mathbf{b_r} \circ \mathbf{c_r}$, $\mathcal{C}$, a given data tensor. The implementation is in Algorithm \ref{lrat}. In a recent work, \cite{fat_nav}, a more practical choice of the parameter $\sigma$ led to a more efficient and accurate algorithm. The practical regularization method is based on a flexible Golub-Kahan (fGK) process. In Figures \ref{fig:my_label4} and \ref{fig:my_label5}, we implemented the LRAT + fGK algorithm with the convergence plot; the sparse model predicts the general phenomena of the Covid-19 cases.

.
\\
In figure $7$, the error plot shows the difference between the original data and estimation. The horizontal line detects an anomaly and this spot coincides with a spike in the number of infections.

\subsection{Hotspot Identification}
 \cite{pmid_hotspot,pmid_hotspot2} Recently, “hotspots” in infectious disease epidemiology have been increasingly used, and they dictate the implementation of appropriate control measures for the specific place. Despite “Hotspots “has not a concrete definition, it is described variously as per area of elevated incidence, prevalence, higher transmission efficiency, or higher chance of disease emergence. Our research is also on Covid-19 pandemic-infected population data, so we defined “hotspots” as the geographical area where the higher intensity of disease prevalence and transmission rate as per population density and flow in the specific area. More specifically, we have defined a threshold in the specific area as per population and its activities which precisely indicates the hotspots there. Identification of hotspots attracts the attention of authorities so that more efficient control measures may be implemented by targeting these areas to sustain further transmission.

\subsubsection{Sparse Optimization for Hotspot Identification}
Our goal is to detect hotspots rapidly. Our goal is to have the following decomposition: $\mathcal{Y} = \mathcal{L} + \mathcal{S}$ where $\mathcal{Y}$ is the given tensor, $\mathcal{L}$ is a low rank reconstructed tensor of $\mathcal{Y}$ and $\mathcal{S}$ is the sparse tensor. In video processing, to detect anomalous activities, the original video is separated into its background and foreground subspaces. The tensor $\mathcal{L}$ is the background and $\mathcal{S}$ is the foreground. The sparse tensor $\mathcal{S}$ can provide anomalous activities. Similarly, $\mathcal{S}$ will contain hotspot occurrence. Thus, we use the sparse tensor model and implement LRAT \ref{lrat}.

\subsubsection{Hotspot Detection with Practical Threshold}
First, we convert the Covid-19 data tensor into a tensor with each entry as the rate of change in the number of infections from the prior week by computing \[ \frac{\text{Number of infections this week} - \text{Last week's total infections} }{\text{Last week's number of infections}}\] for each entry. Then we apply LRAT algorithm to the new tensor. On the newly converted tensor, we apply the practical threshold mean + 5*standard deviation. The resulting plot from LRAT spikes above the threshold line in the week considered to have hotspots. We compare it with the original graph to see that it is consistent, see fig(4.3) below.

\begin{figure}[H]
    \centering
    \includegraphics[scale = 0.6]{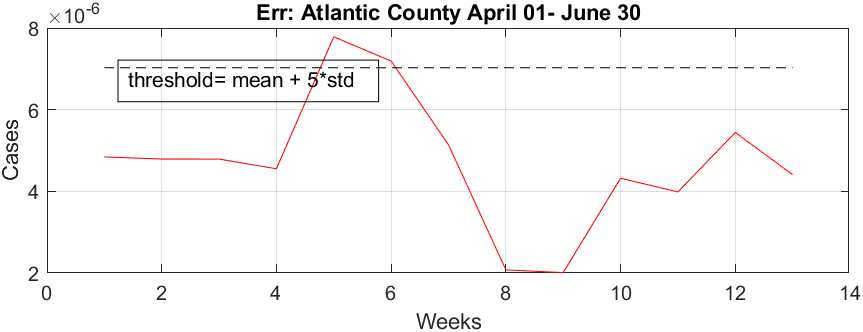}
    \caption{Error Plot for Atlantic County: error between actual vs projected relative data with margin line (mean + 5* standard deviation)}
    \label{fig:my_label6}
\end{figure}

\subsection{Tensor Completion for Predicting Covid-19 Infection Cases}
\begin{figure}[H]
    \centering
    \includegraphics[scale =0.35]{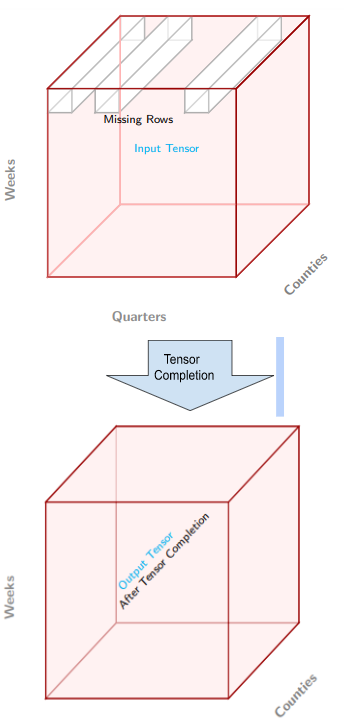}
    \caption{Our Tensor Completion Architecture on Covid-19 data Tensor}
    \label{fig:my_label9}
\end{figure}

The tensor completion is the problem of completing the missing or unobserved entries of the partially observed tensor. Tensor completion algorithms have a wide range of applications in the field of big data \cite{bigdata}, computer vision such as image completion \cite{fat_nav,WangNa1} which focus on filling the missing entries in a presence of noise. Other important applications of tensor completion are link prediction\cite{Liu,Tam} and recommendation system \cite{recom,Taper} and video completion \cite{video}. With the given the tensor $\mathcal{L}$ of order $n$ with missing entries for a given rank, the tensor completion optimization problem can be formulated as the following:
\begin{align*}
    \text{minimize}_{\bf \mathcal{L}} \:\, \text{rank}(\bf \mathcal{L})\\
    \text{subject to}\:\, \mathcal{L} (\Omega) = \mathcal{C}(\Omega) 
\end{align*}
In the work of Wang and Navasca \cite{nav_wang}, this optimization is reformulated to the tensor sparse model. To apply the tensor sparse model with constraints in prediction of Covid-19 infection cases, we set up our data by removing some column data from the original tensor; see Figure \ref{fig:my_label9}. We implement Algorithm \ref{lrat} to complete the Covid-19 tensor with the observed data constraints.  The reconstructed tensor exhibits the same pattern of the tensor cases even though there is some dissimilarities in particular numerical data.

 \begin{algorithm}[H] \label{lrat}
\caption{\cite{nav_wang,fat_nav,fatou}Tensor Completion Via CP with ISTA}
 \begin{algorithmic}
\Require Tensor $\mathcal{C} \in \mathbb{R}^{I \times J \times K}$rank of the tensor R,a regularization parameter $\lambda$ and a scale t$>$0.
\Ensure An approximated tensor $\mathcal{X}$ 
\\
1: \textbf{Initialize} a tensor $\mathcal{C}^0=[\bm{\sigma}^0;\mathbf{A}^0,\mathbf{B}^0,\mathbf{C}^0]_T$.\\
2:\textbf{Update steps:} \\
I) \textbf{Update} factor matrices $\mathbf{A,B,C}$:\\
\hspace*{.5cm}Compute $\mathbf{U}^n$ = $\mathbf{D^n(C\odot B)^T}$ and let $d_n=\max\{\|\mathbf{U}^n{\mathbf{U}^n}^T\|_F,1\}$.\\
\hspace*{.5cm}\textbf{Compute} $\mathbf{D}^n$ and $\mathbf{A}^{n+1}$ by
\begin{align*}
&\mathbf{D}^n=\mathbf{A}^n-\frac{1}{td_n}\nabla_\mathbf{A} f(\mathbf{A}^n,\mathbf{B}^n,\mathbf{C}^n,\bm{\sigma}^n),\\
&\mathbf{A}^{n+1}=\mathbf{D}^n\text{diag}(\|\mathbf{d}_1^n\|,\cdots,\|\mathbf{d}_T^n\|)^{-1}
\end{align*}
\hspace*{.5cm}where $\mathbf{d}_i^n$ is the $i$-th column of $\mathbf{D}^n$ for $i=1,\cdots,T$.\\
\hspace*{.5cm}\textbf{Compute} $\mathbf{V}^n$ as $\mathbf{U^n}$  and let $e_n=\max\{\|\mathbf{V}^n{\mathbf{V}^n}^T\|_F,1\}$.\\
\hspace*{.5cm}\textbf{Compute} $\mathbf{E}^n$ and $\mathbf{B}^{n+1}$ by
\begin{align*}
&\mathbf{E}^n=\mathbf{B}^n-\frac{1}{te_n}\nabla_\mathbf{B} f(\mathbf{A}^{n+1},\mathbf{B}^n,\mathbf{C}^n,\bm{\sigma}^n),\\
&\mathbf{B}^{n+1}=\mathbf{E}^n\text{diag}(\|\mathbf{e}_1^n\|,\cdots,\|\mathbf{e}_T^n\|)^{-1}
\end{align*}
\hspace*{.5cm}where $\mathbf{e}_i^n$ is the $i$-th column of $\mathbf{E}^n$ for $i=1,\cdots,T$.\\
\hspace*{.5cm}\textbf{Compute} $\mathbf{W}^n$ by (1) and let $f_n=\max\{\|\mathbf{W}^n{\mathbf{W}^n}^T\|_F,1\}$.\\
\hspace*{.5cm}\textbf{Compute} $\mathbf{F}^n$ and $\mathbf{C}^{n+1}$ by
\begin{align*}
&\mathbf{F}^n=\mathbf{C}^n-\frac{1}{tf_n}\nabla_\mathbf{C} f(\mathbf{A}^{n+1},\mathbf{B}^{n+1},\mathbf{C}^n,\bm{\sigma}^n),\\
&\mathbf{C}^{n+1}=\mathbf{F}^n\text{diag}(\|\mathbf{f}_1^n\|,\cdots,\|\mathbf{f}_T^n\|)^{-1}
\end{align*}
\hspace*{.5cm}where $\mathbf{f}_i^n$ is the $i$-th column of $\mathbf{F}^n$ for $i=1,\cdots,R$.\\
II) Update the row vector $\bm{\sigma}$:\\
`\hspace*{.5cm}\textbf{Compute} $\mathbf{Q}^{n+1}$ by (2) and let $\eta_n=\max\{\|\mathbf{Q}^{n+1}{\mathbf{Q}^{n+1}}^T\|_F,1\}$.\\
\hspace*{.5cm}\textbf{Compute} $\bm{\beta}^{n+1}$ by (3) and use the soft thresholding:
\hspace*{.5cm} $$\bm{\sigma}^{n+1}=\mathcal{T}_{\frac{\lambda}{t\eta_n}}(\bm{\beta}^{n+1}).$$
3: \textbf{Denote} the limitations by $\mathbf{\hat{A}},\mathbf{\hat{B}},\mathbf{\hat{C}},\bm{\hat{\sigma}}$,
compute $\mathcal{\hat{T}}=[\bm{\hat{\sigma}};\mathbf{\hat{A}},\mathbf{\hat{B}},\mathbf{\hat{C}}]_R$
and count the number $\hat{T}$ of nonzero entries in $\bm{\hat{\sigma}}$.\\
4: \textbf{Impose} constraints $\mathcal{C}(\Omega) = \mathcal{T} (\Omega)$.\\
5: \textbf{Return} The tensor $\mathcal{\hat{T}}$ with the estimated rank $\hat{R}$.
\end{algorithmic}
\end{algorithm}
\vspace{1cm}

\subsubsection{Numerical Experiments on Prediction of Infected Cases}
\label{tc-lrat}
\begin{figure}[H] \label{predict1}
    \centering
    \includegraphics[scale =0.39]{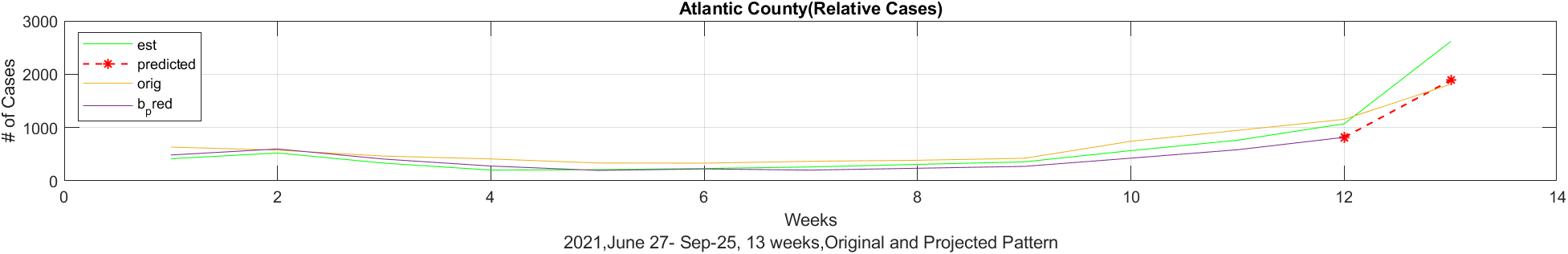}
    \caption{After Tensor Completion Applied to Atlantic County}
\end{figure}

\begin{figure}[H] \label{predict2}
    \centering
    \includegraphics[scale =0.39]{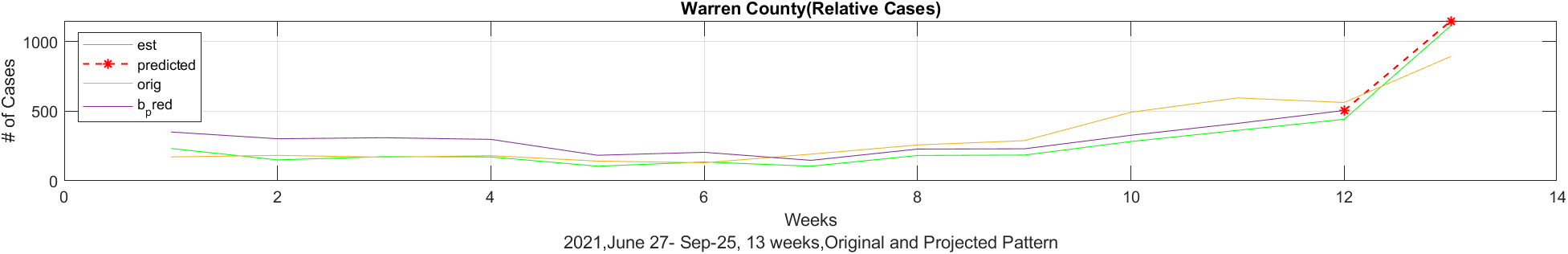}
    \caption{After Tensor Completion Applied to Warren County}
\end{figure}

\begin{figure}[H] \label{predict3}
    \centering
    \includegraphics[scale =0.30]{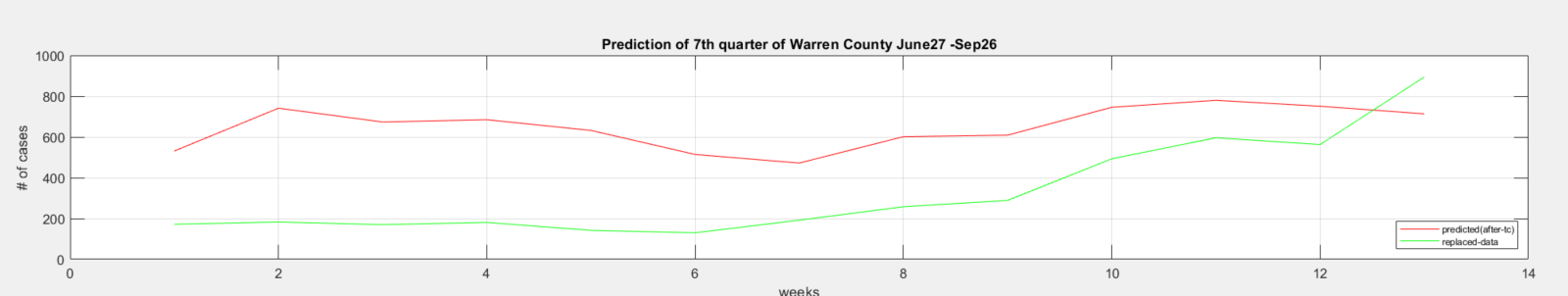}
    \caption{Prediction Covid-19 cases of last quarter based tensor completion formulation: $7^{th}$ quarter of Warren county cases is  replaced by the mean of the remaining data and applied the LRAT \ref{lrat} algorithm.}
\end{figure}

\begin{figure}[H] \label{predict4}
    \centering
    \includegraphics[scale =0.30]{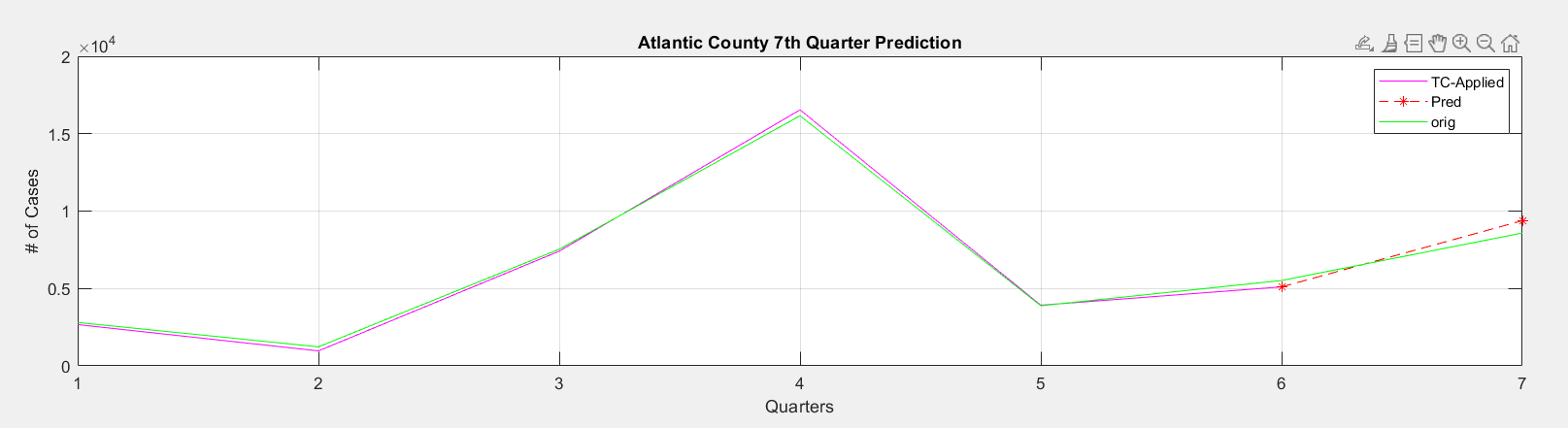}
    \caption{7th Qtr Prediction of Covid Infections: this is the total number of cases prediction of Atlantic county for $7^{th}$ Quarter (June27-Sep25}
\end{figure}

\begin{figure}[H] \label{predict5}
    \centering
    \includegraphics[scale =0.4]{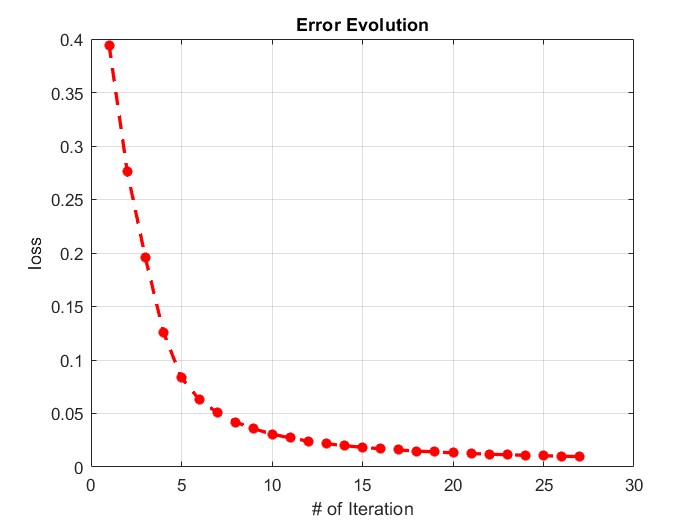}
    \caption{Evolution of error for each epoch}
\end{figure}

We implement a tensor completion algorithm to our tensor data. We replace the last (most current) week's data of the Atlantic and Warren counties with the mean of the remaining data. Then, the missing values on the tensor is completed via low rank approximation of the Covid-19 tensor. We have the following very nice prediction; see Figures \ref{predict1}-\ref{predict5}.






\section{Future Works and Conclusion}
\label{sec:conclusions}
In this work, we apply various tensor models and tensor algorithms to analyze Covid-19 data. The standard tensor models, CP and HOOI, with their off-the-shelves algorithms, ALS and HOOI, are tested against a new sampling method for ALS (SMALS). The numerical results are very promising as it cuts the down time while keeping the Frobenious norm errors relatively consistent with ALS. Here tensor sparse model \cite{nav_wang,fat_nav,} is used as the model for prediction of future Covid-19 infection cases as a tensor completion problem. The numerical results are impressive as the tensor completion algorithm can predict infection a week ahead as well as a quarter ahead. Moreover, the tensor sparse model can locate which counties exhibit hotspots. The sparse tensor model is based on proximal algorithms and flexible hybrid method by Golub-Kahan methods for efficient practical implementation.

In our future work, we would like to have a more mathematical and methodical technique in locating and detecting hotspots. We have used $l_1$ minimization in the tensor completion; we will  work on the efficacy of $l_0$ minimization of low-rank approximation of CP decomposition and tensor completion algorithm in the proximal framework.

\section*{Acknowledgments}
This material is based upon work supported by the National Science Foundation under Grant No. DMS-1439786 while the author, C. Navasca, was in residence at the Institute for Computational and Experimental Research in Mathematics in Providence, RI, during the Model and Dimension Reduction in Uncertain and Dynamic Systems Program. C. Navasca is also in  part supported by National Science Foundation No. MCB-2126374.

\bibliographystyle{siamplain}
\bibliography{tensor_reg}
\end{document}